\documentclass[%
  twocolumn
   , hidempi
]{mpi2015-cscpreprint}

\usepackage[american]{babel}

\usepackage{graphicx}

\usepackage{amsmath}
\usepackage{amssymb}
\usepackage{algorithmic}
\usepackage{algorithm}
\usepackage{verbatim}
\usepackage{xspace,extarrows,fdsymbol}
\usepackage{tikz}
\usetikzlibrary{positioning}

\usepackage{pgfplots}

\pgfplotsset{compat = 1.13,
	colormap name = viridis,
	unbounded coords = jump}

\tikzset{every picture/.style={/utils/exec={\normalfont}}}

\usepackage{caption}
\usepackage{subcaption}
\definecolor{myRed}{HTML}{E34A33}
\definecolor{myBlue}{HTML}{0571B0}
\definecolor{myBrown}{HTML}{A6611A}

\usepackage{xspace,extarrows,fdsymbol}

\usepackage{url}
\usepackage[nameinlink]{cleveref}

\crefformat{equation}{\textup{#2(#1)#3}}
\crefrangeformat{equation}{\textup{#3(#1)#4--#5(#2)#6}}
\crefmultiformat{equation}{\textup{#2(#1)#3}}{ and \textup{#2(#1)#3}}
{, \textup{#2(#1)#3}}{, and \textup{#2(#1)#3}}
\crefrangemultiformat{equation}{\textup{#3(#1)#4--#5(#2)#6}}%
{ and \textup{#3(#1)#4--#5(#2)#6}}{, \textup{#3(#1)#4--#5(#2)#6}}
{, and \textup{#3(#1)#4--#5(#2)#6}}

\usepackage{enumitem}

\usepackage{todonotes}

\usepackage{graphics} % for pdf, bitmapped graphics files

\usepackage{mathtools}
\usepackage{ntheorem}

\usepackage{xspace}

\usepackage[nameinlink]{cleveref}

\newtheorem{theorem}{Theorem}[section]
\newtheorem{problem}{Problem}[section]
\newtheorem{prop}{Proposition}[section]
\crefname{problem}{problem}{problems}
\Crefname{problem}{Problem}{Problems}
\Crefname{prop}{Proposition}{Propositions}
\newtheorem{lemma}{Lemma}[section]
\newtheorem{corollary}[theorem]{Corollary}
\newtheorem{remark}[theorem]{Remark}
\newtheorem{defn}[theorem]{Definition}

\newcommand{\N}{\ensuremath\mathbb{N}}

\newcommand{\R}{\ensuremath\mathbb{R}}
\newcommand{\C}{\ensuremath\mathbb{C}}

\newcommand{\calF}{\mathcal{F}}
\newcommand{\calL}{\mathcal{L}}

\newcommand{\calP}{\mathcal{P}}
\newcommand{\calV}{\mathcal{V}}
\newcommand{\calW}{\mathcal{W}}

\DeclareMathOperator{\Span}{span}
\DeclareMathOperator{\diag}{diag}

\newcommand{\state}{x}
\newcommand{\stateDim}{N}
\newcommand{\inpVar}{u}
\newcommand{\inpVarDim}{m}
\newcommand{\outVar}{y}
\newcommand{\outVarDim}{\ell}

\newcommand{\Ns}{N_s}
\newcommand{\Np}{N_p}

\newcommand{\rhoa}{\rho_a}
\newcommand{\rhow}{\rho_w}
\newcommand{\rhov}{\rho_v}

\newcommand{\param}{p}
\newcommand{\paramCom}[1]{p_{#1}}
\newcommand{\paramSet}{\mathbb{P}}
\newcommand{\paramSetDim}{\nu}
\newcommand{\newparam}{\pi}

\newcommand{\system}{\Sigma}

\newcommand{\red}[1]{\widehat{#1}}
\newcommand{\stateRed}{\red{\state}}

\newcommand{\stateDimRed}{n}
\newcommand{\outVarRed}{\red{\outVar}}
\newcommand{\systemRed}{\red{\Sigma}}
\newcommand{\Ered}{\red{E}}
\newcommand{\Ared}{\red{A}}
\newcommand{\Bred}{\red{B}}

\newcommand{\Cred}{\red{C}}

\newcommand{\outVarLap}{Y}
\newcommand{\inpVarLap}{U}

\newcommand{\transfer}{H}
\newcommand{\transferRed}{\red{\transfer}}
\newcommand{\frequency}{s}

\newcommand{\rank}{\mathsf{rank}}
\newcommand{\pmor}{\textsf{PMOR}\xspace}
\newcommand{\fom}{\textsf{FOM}\xspace}
\newcommand{\mor}{\textsf{MOR}\xspace}
\newcommand{\rom}{\textsf{ROM}\xspace}

\begin{document}
  
%%%%%%%%%%%%%%%%%%%%%%%%%%%%%%%%%%%%%%%%%%%%%%%%%%%%%%%%%%%%%%%%%%%%%%%%%%%%%%%%
% PAPER INFORMATION.                                                           %
%%%%%%%%%%%%%%%%%%%%%%%%%%%%%%%%%%%%%%%%%%%%%%%%%%%%%%%%%%%%%%%%%%%%%%%%%%%%%%%%

\title{Parametric model reduction via rational interpolation along parameters}
  
\author[$\ast$]{Ion Victor Gosea}
\affil[$\ast$]{Max Planck Institute for Dynamics of Complex Technical Systems,
	Sandtorstr. 1, 39106 Magdeburg, Germany.\authorcr
  \email{gosea@mpi-magdeburg.mpg.de}, \orcid{0000-0003-3580-4116}}
  
\author[$\dagger$]{Serkan Gugercin}
\affil[$\dagger$]{Department of Mathematics and Computational Modeling and Data
	Analytics Division, Academy of Integrated Science, Virginia Tech, Blacksburg,
	VA 24061, USA.\authorcr
  \email{gugercin@vt.edu}, \orcid{0000-0003-4564-5999}}
  
\author[$\ddagger$]{Benjamin Unger}
\affil[$\ddagger$]{SC Simulation Technology, University of Stuttgart, Germany.\authorcr
	\email{benjamin.unger@simtech.uni-stuttgart.de}, \orcid{0000-0003-4272-1079}}

\shorttitle{PMOR based on rational interpolation along parameter-dependent curves}
\shortauthor{I. V. Gosea, S. Gugercin, B. Unger}
\shortdate{}
  
\keywords{\small model reduction, parametric systems, rational interpolation, Sylvester equations, holomorphic functions. \normalsize}

%\msc{MSC1, MSC2, MSC3}
  
\abstract{%
    We present a novel projection-based model reduction framework for parametric linear time-invariant systems that allows interpolating the transfer function at a given frequency point along parameter-dependent curves as opposed to the standard approach where transfer function interpolation is achieved for a discrete set of parameter and frequency samples. We accomplish this goal by using parameter-dependent projection spaces. Our main result shows that for holomorphic system matrices, the corresponding interpolatory projection spaces are also holomorphic. The coefficients of the power series representation of the projection spaces can be computed iteratively using standard methods. We illustrate the analysis on three numerical examples.}
  
%     \novelty{This work proposes a new MOR method that .}  
  
\maketitle

%%%%%%%%%%%%%%%%%%%%%%%%%%%%%%%%%%%%%%%%%%%%%%%%%%%%%%%%%%%%%%%%%%%%%%%%%%%%%%%%
% PAPER CONTENT.                                                               %
%%%%%%%%%%%%%%%%%%%%%%

\section{Introduction}

For a parameter vector $\param\in\paramSet\subseteq\R^\paramSetDim$, consider the parametric dynamical system in the state-space form
\begin{equation}
\label{eq:FOM}
\system(\param) \colon \left\{\begin{aligned}
E(\param)\dot{\state}(t;\param) &= A(\param)\state(t;\param) + B(\param)\inpVar(t),\\
\outVar(t;\param) &= C(\param)\state(t;\param),\\
\state(0;\param) &= 0,
\end{aligned}\right.
\end{equation}
with  matrix functions $E,A\colon \paramSet\to\R^{\stateDim\times\stateDim}$, $B\colon \paramSet\to\R^{\stateDim\times\inpVarDim}$, and $C\colon \paramSet\to\R^{\outVarDim\times\stateDim}$. 
We assume that 
$E(\param)$ is nonsingular for every $\param\in\paramSet$.
In \eqref{eq:FOM}, we refer to $\state$, $\inpVar$, and $\outVar$ as the \emph{states}, \emph{inputs}, and \emph{outputs}, respectively.
The parametric dynamical systems of the form \eqref{eq:FOM} arise in many applications ranging from inverse problems to optimal control to uncertainty quantification and 
the parameter vector $\param$ enters the model in various ways, representing, for example, material properties, system geometry, and operating conditions; see, e.g., 
\cite{BGW15surveyMOR,hesthaven2016certified,QuaMN16} and the references therein. 
Our standing assumption is that $\stateDim$ is large and hence simulating \eqref{eq:FOM} for a given input $\inpVar$ and a given parameter $\param$ is expensive. Therefore, 
as required in many prominent applications,
the need to repeat these simulations/computations for many parameter values and input selections  leads to a big computational burden. This is what \emph{parametric model reduction} (\pmor) aims to resolve. The goal of \pmor 
is to replace the  \emph{full-order model} (\fom) \eqref{eq:FOM} by a \emph{reduced-order model} (\rom) of the form
\begin{equation}
\label{eq:ROM}
\systemRed(\param) \colon \left\{\begin{aligned}
\Ered(\param)\dot{\stateRed}(t;\param) &= \Ared(\param)\stateRed(t;\param) + \Bred(\param)\inpVar(t),\\
\outVarRed(t;\param) &= \Cred(\param)\stateRed(t;\param),\\
\stateRed(0;\param) &= 0,
\end{aligned}\right.
\end{equation}
with $\Ered,\Ared\colon\paramSet\to\R^{\stateDimRed\times\stateDimRed}$, $\Bred\colon \paramSet\to\R^{\stateDimRed\times\inpVarDim}$, $\Cred\colon \paramSet\to\R^{\stateDimRed\times\outVarDim}$, and $\stateDimRed\ll\stateDim$ such that the output $\outVarRed(t;\param)$ of the 
\rom approximates the output 
$\outVar(t;\param)$ of the \fom  with high fidelity over a wide range of parameters 
and input selection. 
More precisely, we want the approximation error $\|\outVar-\outVarRed\|$ to be small for any $\inpVar\in\calL_2(0,\infty,\R^{\inpVarDim})$ and any parameter $\param\in\paramSet$. 

\subsection{Projection-based \pmor}
There are plethora of methods to construct the \rom $\systemRed(\param)$, we refer the reader to 
%\cite{BGW15surveyMOR,QuaMN16,hesthaven2016certified,AntBG20,BCOW2017morBook}
\cite{AntBG20,BCOW2017morBook,BGW15surveyMOR,hesthaven2016certified,QuaMN16} for details. 
% \VGcomment{Maybe also cite \cite{HMMS21}...somewhere}
% \SGcomment{What is this citation; it is not showing up}
% \VGcomment{It is the latest MPI work  by Hund et al...see \url{https://arxiv.org/abs/2103.03136}}
% \SGcomment{Victor, I think we should definitely cite that, i was going to do in the future work reference. I kept these citations here only to survey papers and books. But we should cite that work for sure.}
% \BUcomment{I absolutely agree. In particular, since their framework allows to add additional parameter dependencies into the reduced model.}
Common to most of these approaches is that they can be realized via a Petrov-Galerkin framework: Construct two \mor bases  $V,W\in\R^{\stateDim\times\stateDimRed}$ such  that
$\state(t,\param) \approx V \stateRed(t,\param)$. Then, 
substitute this approximation into \eqref{eq:FOM} and enforce a Petrov-Galerkin condition on the residual to obtain
the reduced-order matrices as
\begin{equation}
\label{eq:projection_nonp}
\begin{aligned}
\Ered(\param) &\vcentcolon= W^\top E(\param)V, & \Ared(\param) &\vcentcolon= W^\top A(\param) V,\\
\Bred(\param) &\vcentcolon= W^\top B(\param), & \Cred(\param) &\vcentcolon= C(\param)V.
\end{aligned}
\end{equation}
The task of model reduction is thus essentially equivalent to determining $\stateDimRed$-dimensional subspaces $\calV \vcentcolon= \Span(V)$ and $\calW \vcentcolon= \Span(W)$ of $\R^{\stateDim}$ such that the ROM~\eqref{eq:ROM} obtained via projection onto these spaces is a good approximation of~\eqref{eq:FOM}. Even though it is not the focus of this paper, we note that there are data-driven approaches to 
\pmor in which $\systemRed(\param)$ is constructed without access to the \fom dynamics in \eqref{eq:FOM} and with only access to input-output data; see, e.g., \cite{An12two,CRG20,grimm2018parametric,grivet2017compact,IA14}, and the references therein. 

For the linear parametric dynamical systems \eqref{eq:FOM} and \eqref{eq:ROM} we consider here, the concept of transfer function provides a natural framework to analyze the \mor 
problem. Let $\outVarLap(\frequency,\param)$ and $\inpVarLap(\frequency,\param)$ denote Laplace transforms of $\outVar(t,\param)$ and $\inpVar(t,\param)$. Then,
by taking the Laplace transform of \eqref{eq:FOM}, we obtain
\begin{displaymath}
\outVarLap(\frequency,\param) = \transfer(\frequency;\param) \inpVarLap(\frequency,\param)
\end{displaymath}
where 
\begin{equation} 
\label{eq:transferFunction}
\transfer(\frequency;\param) \vcentcolon= C(\param)\left(\frequency A(\param) - A(\param)\right)^{-1}B(\param)
\end{equation}
is the transfer function of $\system(\param)$.
Similarly, transfer function of the \rom $\systemRed(\param)$ is given by
\begin{equation} 
\label{eq:transferFunctionRed}
\transferRed(\frequency;\param) \vcentcolon= \Cred(\param)\left(\frequency \Ered(\param) - \Ared(\param)\right)^{-1}\Bred(\param).
\end{equation}
In this paper, we will focus on interpolatory approaches to construct 
$\transferRed(\frequency;\param)$. Interpolatory \mor is one of the most commonly employed frameworks to \mor and yield (locally) optimal approximations in the $\mathcal{H}_2$-norm. We skip those details here and refer the reader to \cite{AntBG20}. The interpolatory framework we develop here deviates from the usual approach in the literature as we explain  next.
% Let us emphasize that \Cref{prop:H2optimal} requires that the ROM is an optimal approximation for each parameter. In most of the existing literature on parametric model reduction either the worst-case error over all parameters is minimized (typically via a greedy search) or some kind of least-squares error over the parameter domain is considered. \VGcomment{Add relevant literature here.} 

\subsection{Interpolation problem to construct \texorpdfstring{$\systemRed(\param)$}{the parametric ROM}}

The common approach to interpolatory \pmor chooses $V$ and $W$ so that $\transferRed(s,\param)$ interpolates $\transfer(s,\param)$ at 
some \emph{selected}  right frequency samples $\{\lambda_i\}_{i=1}^{\Ns}$, left frequency samples $\{\mu_i\}_{i=1}^{\Ns}$,
parameter samples $\{\pi_j\}_{j=1}^{\Np}$ along the
right interpolation (tangent) directions $\{r_i\}_{i=1}^{\Ns} \in \C^{\inpVarDim}$ and left interpolation (tangent) directions $\{\ell_i\}_{i=1}^{\Ns} \in \C^{\outVarDim}$; i.e.,
{\small
	$$
	\transfer(\lambda_i,\pi_j)r_i = \transferRed(\lambda_i,\pi_j)r_i~~~\mbox{and}~~~
	\ell_i^\top \transfer(\mu_i,\pi_j) = \ell_i^\top \transferRed(\mu_i,\pi_j)
	$$
}\noindent
for $i=1,2,\ldots,\Ns$ and $j=1,2,\ldots,\Np$. 
One can also enforce interpolating the derivatives of $H$ with respect to $\frequency$ and $\param$, and the discussion here directly extends. However, for brevity, we only focus on  simple interpolation in this paper.
We show in \Cref{thm:rationalInterp} how to construct~$V$ and~$W$ to satisfy the interpolation conditions listed above. These are \emph{discretized} interpolation conditions in the sense that they hold over a discrete set of sampling points.

% In the non-parametric case, i.e., $\paramSet$ consists of a single parameter, \Cref{prop:H2optimal} can be solved by interpolating the transfer function \eqref{eq:transferFunction} at carefully selected frequency points, see for instance~\cite{BeaG17}.
%Thus, as a first step towards solving \Cref{prop:H2optimal} we are faced with the task of interpolating the transfer function along parameter-dependent curves in the frequency domain. More precisely, we are interested in solving the following problem.
In this paper, we consider a more general  problem of \emph{interpolating $\transfer(\frequency,\param)$ along parameter-dependent curves in the frequency domain.}
More precisely, we are interested in solving the following problem.
\begin{problem} \label{prob:newintproblem}
	Consider the dynamical system \eqref{eq:FOM} with transfer function $\transfer(\frequency;\param)$.
	% 	is asymptotically stable for each $\param\in\paramSet$. 
	For given functions $\lambda\colon\paramSet\to\C$, $\mu\colon\paramSet\to\C$,  $r\colon\paramSet\to\C^{\inpVarDim}$,
	and 
	$\ell \colon\paramSet\to\C^{\outVarDim}$,
	construct a \rom %of minimal dimension 
	with transfer function $\transferRed(\frequency;\param)$ that  tangentially interpolates $\transfer$  
	at $\lambda$  along the right tangent directions $r$ and
	at $\mu$ along the right tangent directions 
	$\ell$ for all parameters, i.e., 
	$\transferRed(\frequency;\param)$	 satisfies
	\begin{subequations} \label{eq:interpolationCondition} 
		\begin{align} 
		\transfer(\lambda(\param);\param)r(\param) &= \transferRed(\lambda(\param);\param)r(\param), ~~\mbox{and} \\
		\ell(\param)^\top \transfer(\lambda(\param);\param) &= \ell(\param)^\top\transferRed(\lambda(\param);\param),~~
		\mbox{for~all}~\param\in\paramSet. 
		%\\
		%	\ell(\param)^\top	\transfer'(\lambda(\param);\param)r(\param) &=\ell(\param)^\top \transferRed'(\lambda(\param);\param)r(\param). \nonumber
		\end{align}
	\end{subequations}
\end{problem}

In general, we cannot expect to find constant matrices $V,W\in\R^{\stateDim\times\stateDimRed}$ with small $\stateDimRed$ such that~\eqref{eq:interpolationCondition} is satisfied for all parameters $\param\in\paramSet$. Instead, motivated by the lower-bound for the Kolmogorov $n$-widths \cite[Thm.~3]{UngG19}, we propose to construct parameter dependent model reduction bases $\calV(\param)$ and~$\calW(\param)$, exemplified by the matrix functions
\begin{align*}
V,W\colon\paramSet\to\R^{\stateDim\times\stateDimRed}.
\end{align*}
Our analysis is inspired by the ideas presented in \cite{WitTKAS16}, which studied the balanced truncation method for parametric system.

Once the parameter dependent bases are chosen, the \rom is  constructed via projection onto the spaces given by $\calV(\param)\vcentcolon= \Span(V(\param))$ and $\calW(\param)\vcentcolon= \Span(W(\param))$, i.e.,
\begin{equation}
\label{eq:projection}
\begin{aligned}
\Ered(\param) &\vcentcolon= W(\param)^\top E(\param)V(\param), & \Ared(\param) &\vcentcolon= W(\param)^\top A(\param) V(\param),\\
\Bred(\param) &\vcentcolon= W(\param)^\top B(\param), & \Cred(\param) &\vcentcolon= C(\param)V(\param).
\end{aligned}
\end{equation}

\begin{remark}
	Time- and  state-dependent projection matrices are currently heavily investigated in the efficient approximation of transport-dominated phenomena, where a nonlinear projection framework is used to overcome slowly decaying Kolmogorov $n$-widths, see, e.g.,~\cite{BlaSU20,OhlR13,RimPM2019} and the references therein.
\end{remark}

After this introduction, we recall some preliminary results in \Cref{sec:preliminaries}. Our main contribution is presented in \Cref{sec:rationalInterp} with additional computational details presented in  \Cref{sec:computationalDetails}. 

\paragraph*{Notation} 
Besides standard notation, we use multi-indices, i.e., for $j=(j_1,\ldots,j_{\paramSetDim})\in\N_0^{\paramSetDim}$ and $\param = (\param_1,\ldots,\param_{\paramSetDim})$ we write
\begin{displaymath}
\param^j \vcentcolon= \prod_{i=1}^{\paramSetDim} \param_i^{j_i}.
\end{displaymath}

\section{Preliminaries}
\label{sec:preliminaries}

\subsection{Interpolation conditions}
\label{subsec:interpolationConditions}
Interpolatory model reduction~\cite{AntBG20} constructs reduced-order models whose transfer function interpolates the transfer function of the original model at selected interpolation points. For a fixed parameter $\newparam\in\paramSet$, interpolation via projection can be guaranteed as follows \cite{AntBG20,BaurBBG11}.
\begin{theorem}[Tangential interpolation]
	\label{thm:rationalInterp}
	For a fixed parameter $\newparam\in\paramSet$ consider the \fom~\eqref{eq:FOM} with transfer function $H(s;\newparam)$ and the \rom~\eqref{eq:ROM} with transfer function $\red{H}(s;\newparam)$ constructed as in~\eqref{eq:projection_nonp} using $W,V\in\R^{n\times r}$. For interpolation points $\lambda_0,\mu_0\in\C$, assume that $\lambda_0 E(\newparam)-A(\newparam)$ and $\mu_0 E(\newparam) -A(\newparam)$ are nonsingular. Let $r_0\in\R^m$ and $\ell_0\in\R^m$. 
	\begin{enumerate}
		\item\label{thm:rationalInterp:right} If $(\lambda_0 E(\newparam) -A(\newparam))^{-1}B(\newparam)r_0 \in \Span(V)$, then \\ $\transfer(\lambda_0;\newparam)r_0 = \transferRed(\lambda_0;\newparam)r_0$.
		\item\label{thm:rationalInterp:left} If $(\ell_0^\top C(\newparam)(\mu_0 E(\newparam)-A(\newparam))^{-1})^\top \in \Span(W)$, then $\ell_0^\top \transfer(\mu_0;\newparam) = \ell_0^\top \transferRed(\mu_0;\newparam)$.
		% 		\item If the conditions in~\ref{thm:rationalInterp:right}) and~\ref{thm:rationalInterp:left}) are simultaneously satisfied with $\lambda = \mu$, then $\ell^\top \transfer'(\lambda;\param)r = \ell^\top \transferRed'(\lambda;\param)r$.
	\end{enumerate}
\end{theorem}

It is easy to see (cf.~\cite{GalVV04})  that matrices satisfying the conditions in \Cref{thm:rationalInterp} for driving frequencies $\lambda_i$, $\mu_i$ and tangent directions $r_i$, $\ell_i$ ($i=1,\ldots,\stateDimRed$) can be constructed by solving the two Sylvester equations
\begin{subequations}
	\label{eq:projSylvester}
	\begin{align}
	A(\newparam)V - E(\newparam)V\Lambda + B(\newparam)R &= 0,\\
	W^\top A(\newparam) -M^\top W^\top E(\newparam) + L^\top C(\newparam) &= 0,
	\end{align}
\end{subequations}
for the unknowns $V$ and $W$ where
\begin{subequations}
	\label{eq:projSylvester:matrices}
	\begin{align}
	\Lambda &\vcentcolon= \diag(\lambda_1,\ldots,\lambda_{\stateDimRed}) & M &\vcentcolon= \diag(\mu_1,\ldots,\mu_{\stateDimRed}),\\ 
	R &\vcentcolon= \begin{bmatrix}
	r_1 & \ldots & r_{\stateDimRed}
	\end{bmatrix}, & L &\vcentcolon= \begin{bmatrix}
	\ell_1 & \ldots & \ell_{\stateDimRed}
	\end{bmatrix}.
	\end{align}
\end{subequations}
If the driving frequencies and tangent directions are closed under complex conjugation, then one can use real versions of the matrices in~\eqref{eq:projSylvester:matrices}. For fixed $\newparam\in\paramSet$, the condition guaranteeing the existence and uniqueness of solutions to~\eqref{eq:projSylvester} is well-known, see, e.g.,  \cite[Cha.~6]{Ant05}.

\begin{lemma}
	\label{lem:SylvesterEq}
	For  $\newparam\in\paramSet$, the Sylvester equations~\eqref{eq:projSylvester} have a unique solution if and only if $\lambda_i,\mu_i\not\in \sigma(E(\newparam),A(\newparam))$, where 
	\begin{displaymath}
	\sigma(E(\newparam),A(\newparam)) \vcentcolon = \{s\in \C \mid \rank(sE(\newparam)-A(\newparam)) < \stateDim\}.
	\end{displaymath}
	is the spectrum of the matrix pencil $sE(\newparam)-A(\newparam)$.
\end{lemma}

% \subsection{\texorpdfstring{$\mathcal{H}_2$}{H2}-optimality conditions}
% \BUcomment{Not clear if we need this subsection later on}
% While \Cref{thm:rationalInterp} details the construction of a ROM for given interpolation points $\lambda, \mu$ and tangential directions $r,\ell$, it does not provide a strategy to choose these quantities leading to a high-fidelity ROM. In the context of $\mathcal{H}_2$-optimal reduced-order models, the following theorem provides an implicit definition.

% \begin{theorem}[{$\mathcal{H}_2$ optimality conditions,  {\cite[Thm.~7.7]{BeaG17}}}]
% 	\label{thm:H2optimal}
% 	For fixed $\param\in\paramSet$, let~\eqref{eq:ROM} with semi-simple eigenvalues $\lambda_1,\ldots,\lambda_r$ of\ $\Ered^{-1}(\param)\Ared(\param)$ and transfer function $\transferRed(s;\param) = \sum_{i=1}^r \tfrac{\ell_i r_i^\top}{s-\lambda_i}$ be the best-approximation of~\eqref{eq:FOM} with transfer function $\transfer(s;\param)$ of dimension $\stateDimRed$ with respect to the $\calH_2$-norm.Then,
% 	\begin{subequations}
% 		\label{eqn:H2optimalityConditions}
% 	\begin{align}
% 		\transfer(-\lambda_k;\param)r_k &= \transferRed(-\lambda_k;\param)r_k,\\
% 		\ell^\top \transfer(-\lambda_k;\param) &= \ell^\top \transferRed(-\lambda_k;\param), \text{and}\\
% 		\ell_k^\top \transfer'(-\lambda_k;\param)r_k &= \ell^\top \transferRed'(-\lambda_k;\param)r_k
% 	\end{align} 
% 	\end{subequations}
% 	for $k=1,\ldots,\stateDimRed$.
% \end{theorem}

\subsection{Holomorphic functions}
Our analysis requires that the matrix functions in \eqref{eq:FOM} can be expanded in a power series. If the parameter domain is one-dimensional, this is then equivalent to  the matrix functions being holomorphic (resp. analytic). Since we do not intend to restrict our analysis to a single parameter, we recall the appropriate definitions and results for functions of several parameters. For our presentation we follow~\cite{Mal84} and~\cite{WitTKAS16}.

A function $f\colon\C^\paramSetDim\supseteq \paramSet\to\C$ is called holomorphic in $\param = [\paramCom{j}]\in\paramSet$ if the complex derivative 
\begin{displaymath}
f'(\param) = \lim_{h\to 0} \frac{f(\param + hq) - f(\param)}{h}
\end{displaymath}
exists for any $q\in\C^{\paramSetDim}$. It is said to be holomorphic in $\paramSet$, if it is holomorphic in every $\param\in\paramSet$. Many of the results for the one-dimensional case extend to a higher dimensional domain, such as the Cauchy integral formula. In particular, if $f$ is holomorphic, it can locally be represented via a power series. For the analysis of its domain of convergence, we need the following definition, taken from~\cite{Mal84}.

%The domain of convergence of such a power series in multiple variables is a not a circle any more (as in the case of one variable), but More precisely, there exists vectors $\bs{f}_i\in\mathbb{C}^\paramSetDim$ for $i\in\mathbb{N}_0^\paramSet$ and an open set $\mathbb{U}$ such that
%\begin{displaymath}
%	f(\param) = \sum_{(i_1,\ldots,i_\paramSetDim)\in\mathbb{N}_0^\paramSetDim} \bs{f}_{i_1,\ldots,i_\paramSetDim} \prod_{j=1}^\paramSetDim \paramCom{1}^{i_j}\qquad\text{for each\ }\param\in\mathbb{U}.
%\end{displaymath}
%For the analysis of its domain of convergence, we need the following definition, taken from \cite{Mal84}

\begin{defn}[Reinhardt domain]
	An open set $\Omega\subseteq\C^\paramSetDim$ is called \emph{Reinhardt domain}, if $\param = (\param_1,\ldots,\param_{\paramSetDim})\in\Omega$ implies $(\exp(\imath\theta_1)\param_{1},\ldots,\exp(\imath\theta_\paramSetDim)\param_{\paramSetDim})\in\Omega$ for all $(\theta_1,\ldots,\theta_{\paramSetDim})\in\mathbb{R}^{\paramSetDim}$, where $\imath$ denotes the imaginary unit.
\end{defn} 

\begin{theorem}
	\label{thm:powerSeriesRepresentation}
	Let $\paramSet\subseteq\C^\paramSetDim$ be a connected Reinhardt domain containing $0$ and suppose that $f\colon\paramSet\to\C$ is holomorphic in~$\paramSet$. Then there exist unique $f_i\in\C$ for $i\in\N_0^\paramSetDim$ such that
	\begin{equation}
	\label{eq:powerSeriesRepresentation}
	f(\param) = \sum_{i\in\mathbb{N}_0^\paramSetDim} f_{i} \param^i\qquad\text{for each\ }\param\in\paramSet.
	\end{equation}
\end{theorem}

Note  that for simplicity, we have presented \Cref{thm:powerSeriesRepresentation} solely for the expansion point $\bar{p} = 0$. For practical applications, we may want to use a different expansion point or rescale the parameter domain and the system matrices such that $0$ is included in $\paramSet$.

A question that arises immediately is whether there is an holomorphic version of the implicit mapping theorem available. This is indeed the case. For our analysis, we use the following extension of the implicit mapping theorem~\cite{WitTKAS16}.

\begin{prop}
	\label{prop:implicitMapping}
	Consider a function $\calF\colon\C^{\paramSetDim}\times\C^{n_1\times n_2}\to\C^{n_1\times n_2}$ and suppose there exists $\param_0\in\C^{\paramSet}$ and $X_0\in\C^{n_1\times n_2}$ such that $\calF(\param_0,X_0) = 0$ and $\calF$ is holomorphic around this point. If 
	\begin{displaymath}
	0 = \frac{\partial}{\partial h} \calF(\param_0,X_0 + h D)\big|_{h=0}
	\end{displaymath}
	implies $D=0$, then there exists an neighborhood $\calP\subset\C^{\paramSetDim}$ around $\param_0$ and a holomorphic function $X\colon \calP\to\C^{n_1\times n_2}$ such that
	\begin{displaymath}
	\calF(\param,X(\param)) = 0
	\end{displaymath}
	for all $\param\in\paramSet$.
\end{prop}

\section{Rational interpolation along parameter-dependent curves}
\label{sec:rationalInterp}

In this section, we establish the main result that guarantees existence of holomorphic functions $V(p)$ and $W(p)$ such that the reduced model in \eqref{eq:projection} solves the new interpolation problem defined in \Cref{prob:newintproblem}. 
\begin{theorem}
	\label{thm:rationalInterpParametric}
	Consider the dynamical system~\eqref{eq:FOM} and assume that $E,A,B,C$ are holomorphic in the compact set $\paramSet\subseteq\C^{\paramSetDim}$. Assume that for $i=1,\ldots,\stateDimRed$ the holomorphic functions $\lambda_i,\mu_i\colon\paramSet\to\C$ are such that
	\begin{displaymath}
	\lambda_i(\param),\mu_i(\param)\not\in\sigma(E(\param),A(\param))
	\end{displaymath}
	for all $p\in\paramSet$. If the tangent directions $r_i\colon\paramSet\to\C^{\inpVarDim}$ and $\ell_i\colon\paramSet\to\C^{\outVarDim}$ are holomorphic, then there exists holomorphic functions $V,W\colon\paramSet\to\C^{\stateDim\times\stateDimRed}$ satisfying
	\begin{align}
	\label{eq:Sylvester:param:V}A(\param)V(\param) - E(\param)V(\param)\Lambda(\param) + B(\param)R(\param) &= 0,\\
	\label{eq:Sylvester:param:W}W^\top\!(\param) A(\param) -M^\top\!(\param) W^\top\!(\param) E(\param) + L^\top\!(\param) C(\param) &= 0
	\end{align}
	for all $p\in\paramSet$, where $\Lambda(\param)$, $M(\param)$, $R(\param)$, $L(\param)$ are defined as in~\eqref{eq:projSylvester:matrices}, but now with parametric dependence.
	% 	\SGcomment{The equation~\eqref{eq:projSylvester:matrices} does not have parametric dependencies. But it is probably clear.} \BUcomment{I think it should be clear enough. But if you want to add the equations with the parameter-dependency, then it is also okay with me.}
\end{theorem}

{\textit{Proof:}
	We show the assertion only for $V$. The proof for $W$ follows similarly. Define the holomorphic function
	\begin{gather*}
	\calF\colon \paramSet\times \C^{\stateDim\times \stateDimRed}\to \C^{\stateDim\times\stateDimRed},\\
	(\param,V) \mapsto A(\param)V - E(\param)V\Lambda(\param) + B(\param)R(\param).
	\end{gather*}
	Let $p_0\in\paramSet$. Then,  using \Cref{lem:SylvesterEq}, there exists $V_0\in \C^{\stateDim\times\stateDimRed}$ satisfying the condition $\calF(\param_0,V_0) = 0$. In addition, for any $\widetilde{V}\in\C^{\stateDim\times\stateDimRed}$ we obtain
	\begin{align*}
	\frac{\partial}{\partial \varepsilon}\calF(\param,V_0+\varepsilon\widetilde{V}) = A(\param)\widetilde{V} - E(\param)\widetilde{V}\Lambda(\param).
	\end{align*}
	From \Cref{lem:SylvesterEq} we conclude that $\tfrac{\partial}{\partial \varepsilon}\calF(\param,V_0+\varepsilon\widetilde{V}) = 0$ if and only if $\widetilde{V} = 0$. Thus, \Cref{prop:implicitMapping} implies that there exists a neighborhood $\calP\subseteq\C^{\paramSetDim}$ around $\param_0$ and a holomorphic function $V\colon\calP\cap\paramSet\to\C^{\stateDim\times\stateDimRed}$ satisfying $\calF(\param,V(\param)) = 0$. Let $\calP$ denote the maximal neighborhood such that the previous construction holds. It remains to show that $\calP \cap \paramSet = \paramSet$.  Assume $\calP \cap\paramSet \neq \paramSet$ and let $\pi\in \paramSet \setminus \calP$. Repeating the construction, we obtain a neighborhood $\tilde{\calP}\subseteq \C^{\paramSetDim}$ and holomorphic function $\tilde{V}\colon\tilde{\calP}\cap\paramSet\to\C^{\stateDim\times\stateDimRed}$ satisfying $\calF(\param,\tilde{V}(\param)) = 0$. Assume first $\calP\cap\tilde{\calP}\neq\emptyset$. Then there exists $\tilde{\pi}\in \calP\cap\tilde{\calP}$. Due to \Cref{lem:SylvesterEq} and the assumptions we conclude $V(\tilde{\pi}) = \tilde{V}(\tilde{\pi})$. From the holomorphic identity theorem \cite[Thm.~1.2.14]{Nog16} we infer $V = \tilde{V}$, a contradiction. If, on the other hand, $\calP\cap\tilde{\calP}=\emptyset$, we can select further points in $\paramSet$ until we obtain an open covering of $\paramSet$. Since $\paramSet$ is compact, we can choose a finite covering and proceed as before. We conclude $\calP \cap \paramSet = \paramSet$.}
%\end{prof}

%\begin{corollary}
%	If $V$ and $W$ from \Cref{thm:rationalInterpParametric} have constant rank $r$ for all $\param\in\paramSet$, then there exists holomorphic functions $\widetilde{V},\widetilde{W}\colon\paramSet\to\C^{\stateDim\times r}$ with constant rank $r$ for all $\param\in\paramSet$, $\Span(V) = \Span(\widetilde{V})$, $\Span(W) = \Span(\widetilde{W})$, and $\widetilde{V}(\param)$ and $\widetilde{W}(\param)$ are orthogonal for each $\param\in\paramSet$. \BUcomment{Can we do this such that $W^\top\!(\param)E(\param)V(\param) = I_r$ for all $p\in\paramSet$?}
%\end{corollary}
%
%\begin{proof}
%	This is a simple consequence of the smooth SVD \cite{BunBMN91} or smooth QR-decomposition \cite{Rhe88}.
%	\todo{Check details and what we need to extend that to higher-dimensional parameter space}
%\end{proof}

\begin{corollary}
	Suppose that the assumptions enforced in \Cref{thm:rationalInterpParametric} are satisfied and construct a \rom as in~\eqref{eq:projection}. Then the \rom satisfies the interpolation conditions~\eqref{eq:interpolationCondition} for all $\param\in\paramSet$, thus solving \Cref{prob:newintproblem}. 
\end{corollary}

\begin{remark}
	Using \cite[Prop.~3.24]{SchUBG18}, \Cref{thm:rationalInterpParametric} can be extended to structured systems with a transfer function of the form $\transfer(s;\param) = C(\param)(\sum_{k=1}^K h_k(s;\param)A_k(\param))^{-1}B(\param)$, which includes, for instance, delay equations, fractional systems, and viscoelastic dynamics.
\end{remark}

\section{Computational details}
\label{sec:computationalDetails}

Even though we have established the theoretical framework for constructing $V(p)$ and $W(p)$ to solve the new parametric interpolation problem, for a numerically efficient \pmor framework we need to consider the computational aspects in solving \eqref{eq:Sylvester:param:V} and \eqref{eq:Sylvester:param:W}, and performing the projection \eqref{eq:projection}. 

For the brevity of presentation we restrict ourselves in this section to standard state-space systems with $E(\param)\equiv I_{\stateDim}$. For a parameter-dependent $E$ matrix, the construction is similar, but the formulas are more involved.

\subsection{Numerical construction of \texorpdfstring{$V(p)$ and $W(p)$}{the projection matrices}}
\label{subsec:numericalConstructionOfVW}

Assuming that holomorphic matrix functions, \Cref{thm:powerSeriesRepresentation} ensures that we can decompose these matrices as
\begin{gather*} 
A(\param) = \sum_{i\in \N_0^\paramSetDim} \param^i A_i, ~\ 
\Lambda(\param) = \sum_{i\in\N_0^\paramSetDim} \param^i \Lambda_i, ~\ 
V(\param) = \sum_{i\in \N_0^\paramSetDim} \param^i V_i\\
B(\param) = \sum_{i\in \N_0^\paramSetDim} \param^i B_i,\quad\text{and}\quad 
R(\param) = \sum_{i\in \N_0^\paramSetDim} \param^i R_i. 
\end{gather*}

In many practical applications, the system matrices are directly available in such a form with a finite number of terms.
Then the Sylvester equation~\eqref{eq:Sylvester:param:V} becomes
\begin{align*}
0 &= \sum_{j\in \N_0^\paramSetDim} \sum_{i\in \N_0^\paramSetDim} (A_iV_j - V_j\Lambda_i + B_iR_j)\param^{i+j}\\
&= \sum_{\rho\in \N_0^\paramSetDim} \sum_{\substack{i+j=\rho\\i,j\in\N_0^\paramSetDim}} (A_iV_j - V_j\Lambda_i + B_iR_j)\param^{\rho}.
\end{align*}
Using the holomorphic identity theorem \cite[Thm.~1.2.14]{Nog16}, we conclude that for $\rho\in\N_0^\paramSetDim$ we have
\begin{equation}
\label{eq:Sylvester:param:V:coefficient}
\begin{aligned}
0 &= \sum_{\substack{i+j=\rho\\i,j\in\N_0^\paramSetDim}} (A_iV_j - V_j\Lambda_i + B_iR_j)\\
&= A_0V_\rho - V_\rho\Lambda_0 + \sum_{\substack{i+j=\rho\\j\neq\rho}} A_iV_j - V_j\Lambda_i + \sum_{i+j=\rho} B_iR_j,
\end{aligned}
\end{equation}
which provides an iterative method to solve for the coefficients $V_i$. A similar strategy can be obtained for the coefficients for $W$, which we omit here to avoid redundancy.
\begin{corollary}
	Under the assumptions of \Cref{thm:rationalInterpParametric} the Sylvester equation~\eqref{eq:Sylvester:param:V:coefficient} is uniquely solvable for each $\rho\in\N_0^\paramSetDim$.
\end{corollary}
{\hspace{10mm} \textit{Proof:}
	This follows immediately from $A_0 = A(0)$, $\Lambda_0 = \Lambda(0)$, and \Cref{lem:SylvesterEq}. 
}

\vspace{2mm}

Note that if the coefficients of $A$, $\Lambda$, $B$, and $R$ are real (i.e., the interpolation frequencies and tangent directions are closed under conjugation), then the $V_j$ are real thus yielding a real-valued matrix $V(\param)$ for each real parameter $\param\in\paramSet$. 

In numerical computations, we cannot compute all the coefficients $V_i$ and thus have to truncate the power-series expansion at an index based on a tolerance. In other words, for a given tolerance $\tau$, we truncate the power series expansion when $\max_{\param\in\paramSet}|\param^i|\| V_i \|\leq \tau$, and similarly for $W(p)$.  As a consequence, we cannot ensure exact interpolation any longer. A similar issue arises in the usual interpolatory model reduction framework when the required subspace vectors in \Cref{thm:rationalInterp}, namely   
\small
$$
(\lambda_0 E(\newparam) -A(\newparam))^{-1}B(\newparam)r_0, \ \text{and} \ (\ell_0^\top C(\newparam)(\mu_0 E(\newparam)-A(\newparam))^{-1})^\top,
$$
\normalsize
 are computed via iterative solves; see, e.g., \cite{beattie2012inexact}. We revisit this issue in 
\Cref{sec:conc}.

\subsection{Constructing the reduced matrices}
\label{sec:ConstrRedMat}

For simplicity, we only focus on
$\Ared(p)$ in \eqref{eq:projection}; but the discussion extends directly to other reduced order quantities. 

We will work with the truncated quantities, i.e.,
\begin{equation} 
\label{eq:truncatedform}
\begin{gathered}
A(\param) = \sum_{\|k\|\leq \rhoa} \param^k A_k, \\
W(\param) = \sum_{\|j\|\leq \rhow} \param^i W_i, \qquad
V(\param) = \sum_{\|i\|\leq \rhov} \param^i V_i.
\end{gathered}
\end{equation}
For every new parameter vector $\newparam \in \paramSet$,
forming $V(\newparam)$ (and $W(\newparam)$) can be efficiently done using the truncated form as in
\eqref{eq:truncatedform}. However constructing
$\Ared(\newparam)$ requires computing
$
\Ared(\newparam) = W(\newparam)^\top
A(\newparam) V(\newparam),
$
which involves two matrix multiplications in the original dimension $\stateDim$. We resolve this issue 
using the truncated forms
\eqref{eq:truncatedform}:
\begin{equation} \label{eq:termijk}
\Ared(\newparam) 
= \sum_{\|j\|\leq \rhow}\,\sum_{\|k\|\leq \rhoa}\,\sum_{\|i\|\leq \rhov} \left(W_j^\top A_k V_i \right) \newparam^{i+j+k}.
\end{equation}
Note that the reduced coefficients  $W_j^\top A_k V_i \in \R^{\stateDimRed}$ in \eqref{eq:termijk} can be precomputed (in the offline stage).
Assuming $\rhoa,\rhov$ and $\rhow$ are modest integers, storing all the coefficients $W_j^\top A_k V_i$ and then forming  the overall sum
can be efficiently computed in the online stage.

\section{Numerical examples}
We illustrate the theoretical analysis on three models.

\subsection{A toy example}
% \BUcomment{Since I removed the toy example from the introduction, what about removing this example here as well? Two examples should be sufficient.}
% \VGcomment{I suggest that we keep it here for the first submission since everything is done explicitly (in case the reviewers do not understand what we are proposing) - we could remove it for the final submission, provided that the paper actually get's accepted...}

Consider a simple example for which the dimension of the parameter set is $\paramSetDim=1$ (and the parameter enters only in the vector $B$). The matrices are as follows
\begin{equation}
\label{eq:toy_mat_1}
\begin{aligned}
A(\param) &\equiv -\text{diag}(1,1,2), & C(\param) &\equiv \begin{bmatrix} 2 & 1 & 1 \end{bmatrix}, \\
B(\param) &= \begin{bmatrix} p & 1-p & 1 \end{bmatrix}^\top, & E(\param) &\equiv I_3.
\end{aligned}
\end{equation}
Hence, it follows that
\begin{align*}
A_0 &= -\text{diag}(1,1,2), \ C_1 = \left[ \begin{matrix} 1 & 1 & 2 \end{matrix} \right], \ A_i, C_i = 0 \ \forall i \geq 1 \\
B_0 &= \left[ \begin{matrix} 0 & 1 & 1 \end{matrix} \right]^\top, \ \    B_1 = \left[ \begin{matrix} 1 & -1 & 1 \end{matrix} \right]^\top, \ B_i = 0 \ \forall i \geq 2.
\end{align*}
Originally, note that $N=3$ and that we choose $n=2$ as the reduction order. Choose interpolation points and tangent directions that are independent of the parameter, e.g.,
\begin{equation}
\lambda_1 = 1, \ \lambda_2 = 3, \ \mu_1 = 2, \ \mu_2 = 4,
\end{equation}
and also $r = \ell =  [1;1]$. Choose the following:
\begin{align*}
\Lambda_0 = \left[ \begin{matrix} 1 & 0 \\ 0 & 3 \end{matrix} \right], \        M_0 = \left[ \begin{matrix} 2 & 0 \\ 0 & 4 \end{matrix} \right], \  \Lambda_i = 0,  \ M_i = 0 \ \forall i \geq 1 \\
R_0 = L_0^\top = \left[ \begin{matrix} 1 & 1 \end{matrix} \right], \ R_i = L_i^\top = 0, \forall i \geq 1.
\end{align*}
For $\rho = 0$, it follows that the equation (\ref{eq:Sylvester:param:V:coefficient}) simplifies to $A_0 V_0-V_0 \Lambda_0+B_0 R_0 = 0$. Similarly, based also on (\ref{eq:Sylvester:param:V:coefficient}), $V_1$ satisfies the following Sylvester equation
\begin{align}
A_0 V_1-V_1 \Lambda_0+ A_1 V_0-V_0 \Lambda_1+B_1 R_0+B_0 R_1 = 0,
\end{align}
which simplifies to $A_0 V_1-V_1 \Lambda_0+B_1 R_0 = 0$. Hence, explicitly compute the first two Taylor coefficients
\begin{align}\label{toy_V0_V1}
\begin{split}
V_0 &= \left[ \begin{matrix} 0 & 0 \\[2mm] \frac{1}{2} & \frac{1}{4} \\[2mm] \frac{1}{3} & \frac{1}{5} \end{matrix} \right], \    V_1 = \left[ \begin{matrix} \frac{1}{2} & \frac{1}{4} \\[2mm] -\frac{1}{2} & -\frac{1}{4} \\[2mm] 0 & 0 \end{matrix} \right],
\end{split}
\end{align}
and $V_i = 0,  \forall i \geq 2$.
%Note that, in this simplified case, one could directly compute $V(\param)$ in (\ref{toy_V0_V1}) as follows
%\begin{align*}
%    V(\param) 
%    &= \left[  \begin{matrix} (\lambda_1 E - A_0)^{-1} B(\param) &  (\lambda_2 E - A_0)^{-1} B(\param) \end{matrix} \right]
%\end{align*}
Next, compute matrix $W^\top = W_0^\top \in \C^{2 \times 3}$ by solving $W_0^T A_0-M_0 W_0^T +L_0 C_0 = 0$, as
\begin{align}
W^\top =   \left[ \begin{array}{ccc} \frac{2}{3} & \frac{1}{3} & \frac{1}{4}\\[2mm] \frac{2}{5} & \frac{1}{5} & \frac{1}{6} \end{array}\right],
\end{align}
and put together the following reduced realization that does indeed depend on the parameter $p$ as follows
\begin{align}\label{eq:realiz1}
\hat{E}^{(1)}(\param) &=   W^\top E V(\param) =\left[\begin{array}{cc} \frac{p}{6}+\frac{1}{4} & \frac{p}{12}+\frac{2}{15}\\[1mm] \frac{p}{10}+\frac{7}{45} & \frac{p}{20}+\frac{1}{12} \end{array}\right], \nonumber \\[1mm]
\hat{A}^{(1)}(\param) &=   W^\top A(\param) V(\param) = \left[ 
\begin{array}{cc} -\frac{p}{6}-\frac{1}{3} & -\frac{p}{12}-\frac{11}{60}\\[1mm] -\frac{p}{10}-\frac{19}{90} & -\frac{p}{20}-\frac{7}{60} \end{array}
\right], \nonumber \\     \hat{B}^{(1)}(\param) &=   W^\top B(\param) = \left[\begin{array}{c} \frac{p}{3}+\frac{7}{12}\\[1mm] \frac{p}{5}+\frac{11}{30}  \end{array}\right], \\  \hat{C}^{(1)}(\param) &=  C(\param) V(\param) = \left[\begin{array}{cc}  \frac{p}{2}+\frac{5}{6} & \frac{p}{4}+\frac{9}{20} \end{array}\right]. \nonumber
\vspace{-2mm}
\end{align}
%Note that, for $p\neq -1$, the realization in (\ref{realiz1}) is equivalent to another for which the parameter enters only in the vector $C$.
%\begin{align}
%    \begin{split}
%           \hat{E}^{(2)}(\param) &= \text{diag}(1,1), \\
%            \hat{A}^{(2)}(\param) &=   \left[\hat{E}^{(1)}(\param)\right]^{-1}  \hat{A}^{(1)}(\param) = \left[\begin{array}{cc} 4 & 3\\ -10 & -7 \end{array}\right], \\
%              \hat{B}^{(2)}(\param) &=   \left[\hat{E}^{(1)}(\param)\right]^{-1}  \hat{B}^{(1)}(\param) = \left[\begin{array}{cc} - 3\\ 10 \end{array}\right], \\
%           \hat{C}^{(2)}(\param) &=    %\hat{C}^{(1)}(\param) = \left[\begin{array}{cc}  \frac{p}{2}+\frac{5}{6} & \frac{p}{4}+\frac{9}{20} \end{array}\right]. 
%    \end{split}
%\end{align}

{We note that the system in (\ref{eq:realiz1}) interpolates the original one in (\ref{eq:toy_mat_1}) at the selected frequencies for \emph{every} value of the parameter $p$.
	We also note that the system in (\ref{eq:realiz1}) is equivalent to a minimal realization of (\ref{eq:toy_mat_1}) for $p \in \{0,1\}$.}

%\SGcomment{Shouldn't we also add that we interpolate for every $\param$ at the selected interpolation points. }

\subsection{Another toy example}

\begin{figure}[ht]
	\hspace{-4mm}	\includegraphics[width=1.1\linewidth]{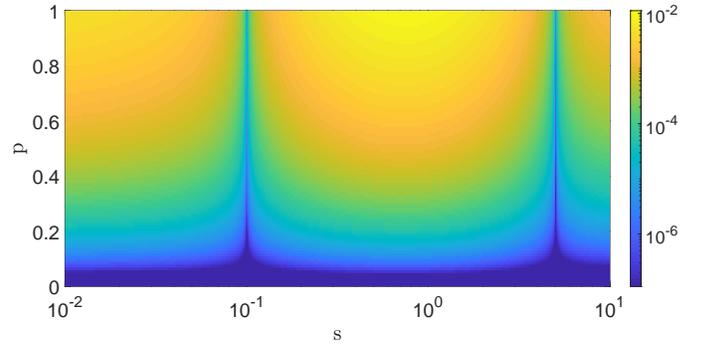}
	\vspace{-6mm}
	\caption{Approximation errors on a 2D grid $(s,p)$.}
	\label{fig:1}
	\vspace{-3mm}
\end{figure}

Consider the following example:
\begin{align}\label{eq:toy_mat2}
\begin{split}
A(\param) &= \left[\begin{array}{ccc} -2 & p & 0\\ -p & -1 & 0\\ 0 & 0 & -1 \end{array}\right], \ \ C(\param)= \left[ \begin{matrix} 1 & 0 & 1 \end{matrix} \right], \\
B(\param) &= \left[ \begin{matrix} 1 & 0 & 1 \end{matrix} \right]^\top, \ \ E(\param) = I_3.
\end{split}
\end{align}
%\BUcomment{Maybe we should change the (1,1) or (3,3) entry of $A(p)$ to something different, since with the current setting, the controllability matrix has rank 1, and thus even after truncation, the system is not minimal. If I'm not mistaken, then for $p=0$ the pencil $(\Ered(p),\Ared(p))$ is singular.}
%\VGcomment{Did that!}
Hence, it follows that:
\begin{align*}
A_0 = \begin{bmatrix} -2 & 0 & 0\\ 0 & -1 & 0\\ 0 & 0 & -1 \end{bmatrix}, \ \ A_1 = \begin{bmatrix} 0 & 1 & 0\\ -1 & 0 & 0\\ 0 & 0 & 0 \end{bmatrix},
\end{align*}
and $A_i=0$ for all $i \geq 2$. Then, we have also that
\begin{equation*}
B_0^\top = C_0 = \begin{bmatrix} 1 & 0 & 1 \end{bmatrix}, \ \text{and} \    B_i^\top = C_i = 0, \forall i \geq 1. 
\end{equation*}
For this case, consider two right interpolation points as:
\begin{align*}
\Lambda(\param) = \begin{bmatrix} 0.1 & 0 \\ 0 & 5 \end{bmatrix} = \Lambda_0, \ \text{and} \ \Lambda_i = 0, \ \forall i \geq 1.
\end{align*}
Note that in this case we use $W^\top(\param) = V^\top(\param)$. The right directions are all ones and the Sylvester equations in (\ref{eq:Sylvester:param:V:coefficient}) simplify to the following collection:
\begin{align}\label{eq:Sylvester_tot2}
\begin{split}
& A_0 V_0-V_0 \Lambda_0+B_0 R_0 = 0, \\
& A_0 V_i-V_i \Lambda_0+A_1 V_{i-1} = 0, \ \forall i \geq 1.
\end{split}
\end{align}
Hence, one can iteratively compute $V_i$ for any positive value of $i$. We do that for all values of $i$ until $\Vert V_i \Vert < \tau $, for a tolerance value of $\tau = 10^{-5}$. This corresponds to a number of 26 Taylor coefficients that need to be computed. Finally, as described in \Cref{sec:ConstrRedMat}, we put together the reduced-order matrices and evaluate the approximation errors for a 2D grid consisting in values $p \in [0,1]$, and $s \in [10^{-2},10^1]$. The results are presented in \Cref{fig:1}. 

%Hence, approximate $V(\param)$ by $\tilde{V}(\param)$, where $\tilde{V}(\param) = \sum_{i=0}^{25} p^i V_i$. 

\subsection{A more involved numerical example}

We analyze the dynamical system originally proposed in \cite{Penzl06} and later modified in \cite{IA14,CRG20} to add a parameter dependence. The dynamics are characterized by the following equations:
\begin{equation}
\label{eq:PenzlEx}
\system(\param) \colon \left\{\begin{aligned}
\dot{\state}(t;\param) &= A(\param)\state(t;\param) + B(\param)\inpVar(t),\\
\outVar(t;\param) &= C(\param)\state(t;\param),\\
\state(0;\param) &= 0
\end{aligned}\right.
\end{equation}
where  $\param \in \paramSet = [0,1]$
and
$A\colon \paramSet\to\R^{1006\times 1006}$
\begin{align*}
A(\param) &= \text{diag} \left( T_1(\param), T_2, T_3, T_4 \right), \ \text{with} \\
T_1(\param) &= \left[ \begin{matrix} -1 & p+100 \\ -100-p & -1  \end{matrix} \right],   T_2 = \left[ \begin{matrix} -1 & 200 \\ -200 & -1  \end{matrix} \right], \\
T_3 &= \left[ \begin{matrix} -1 & 400 \\ -400 & -1  \end{matrix} \right], \ \ T_4 = -\text{diag} \left(1,2,\ldots,1000 \right).
\end{align*}
Additionally, the constant vectors $B$ and $C$ are given by
\begin{equation*}
B = C^T = [10e_6 ; e_{1000}],
\end{equation*}
% \begin{equation}
%     B = C^\top = \begin{bmatrix}1\\0\end{bmatrix} \otimes e_{503},
% \end{equation}
where $e_k$ denotes the $k$-dimensional vector of ones.
% and $\otimes$ the Kronecker product.

Next, we choose 40 logarithmically-spaced interpolation points $\lambda_1,\ldots,\lambda_{40}$ in the interval $[10^{-1},10^3]\imath$ (we are using a one-sided interpolation scheme). Additionally, let the tolerance value be $\tau = 10^{-7}$.
It follows that we need to compute the first $11$ Taylor coefficients of $V(\param)$, i.e., $V_1,V_2,\ldots, V_{10}$, since $V_{11}<\tau$. As in the previous example, use $W(\param) = V\top(\param)$ as left projection matrix, and follow the formulas presented in \Cref{sec:ConstrRedMat}, to compute the corresponding reduced-order matrices.

% \SGcomment{Is this one sided or two-sided? Do we have $\mu_i$'s here?}

\begin{figure}[ht]
	\hspace{-4mm}		
	\includegraphics[width=1.1\linewidth]{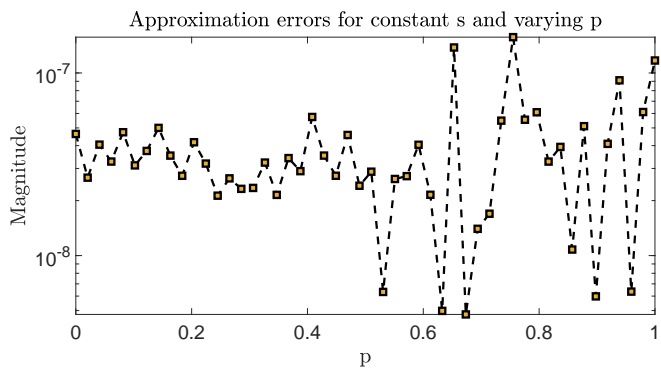}
	\vspace{-6mm}
	\caption{Approximation errors for $s=\lambda_{20}$ and varying $p$ in $[0,1]$.}
	\label{fig:5}
	\vspace{-3mm}
\end{figure}

First we fix the frequency parameter as $s = \lambda_{20} = 8.8862\imath$ and vary $p$ in between $0$ and $1$ (50 linearly-spaced points). We depict the approximation errors for different values of $p$ in \Cref{fig:5}. 
We note that the interpolation errors due to the truncation of the power series are small, in the interval $(10^{-9},10^{-6})$, in accordance with the tolerance $\tau = 10^{-7}$.

%\SGcomment{And this graph shows that ...} 

\begin{figure}[ht]
	\hspace{-4mm}		
	\includegraphics[width=1.1\linewidth]{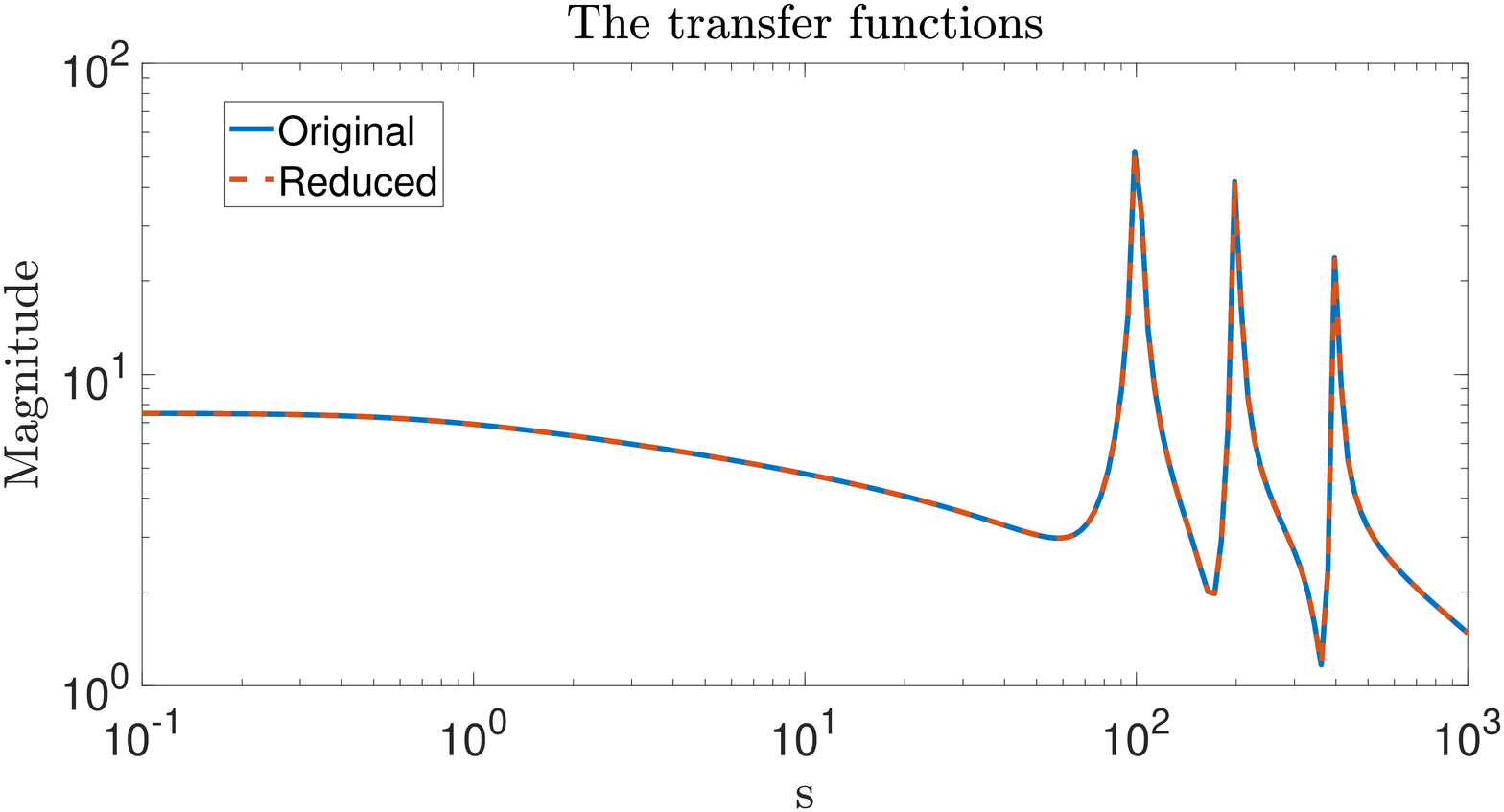}
	\vspace{-6mm}
	\caption{The two transfer functions for $p=0.5$ and $s$ in $[10^{-1},10^3]\imath$.}
	\label{fig:6}
	\vspace{-3mm}
\end{figure}

For the next experiment, we fix the $p$ parameter, i.e., choose $p = 0.5$ and vary the frequency parameter $s$ in the interval $[10^{-1},10^3]\imath$ (200 logarithmically-spaced points). We depict the magnitudes of the two transfer functions (original and reduced) evaluated for different values of $s$ in \Cref{fig:6}, illustrating that \fom  response is indeed well matched.

%\SGcomment{And this figure illustrates that ...} 

%Now, choose  $T_1(\param) = \left[ \begin{matrix} -1 & p+10 \\ -10-p & -1  \end{matrix} \right]$ and keep $p=0.5$ - see the results depicted in Fig.\;\ref{fig:7}.

%\begin{figure}[ht]
%	\begin{center}			
%	\includegraphics[width=\linewidth]{figures/plot_Penzl4.eps}
%	\end{center}
%	\vspace{-3mm}
%	\caption{Approximation errors for $p=0.5$ and varying $s$ in $[10^{-1},10^3]$.}
%	\label{fig:7}
%	\vspace{-3mm}
%\end{figure}

\begin{figure}[ht]
	\hspace{-4mm}			
	\includegraphics[width=1.1\linewidth]{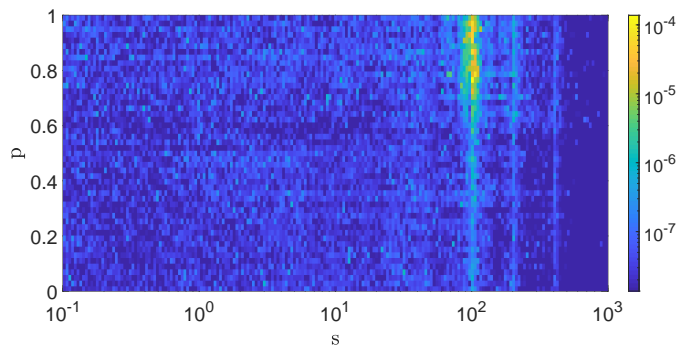}
	\vspace{-6mm}
	\caption{Approximation errors on a 2D grid $(s,p)$.}
	\label{fig:8}
	\vspace{-3mm}
\end{figure}

Finally, we construct a 2D grid consisting in pairs of parameters $(s,p)$ evaluated on the Cartesian product of the two previously-mentioned discrete sets. Then, for all the $200 \times 50 = 10^4$ pairs, we compute the approximation error. The results are presented in \Cref{fig:8}.

%\BUcomment{We should probably choose the tolerance, such that $\max_{\param\in\paramSet} \|\param^i\| \|V_i\| \leq \tau$, because this somewhat reflects the worst truncation error.}

\section{Conclusions and future work}
\label{sec:conc}

We have presented a theoretical framework that allows to construct a \rom whose transfer function interpolates the transfer function of the original high-dimensional system at parameter-dependent interpolation frequencies
along some parameter-dependent directions. The associated parametric projection spaces are proven to have a holomorphic dependency on the parameter and the coefficients of its power series can be computed iteratively using standard methods. 

There are many natural avenues to investigate further. For example, interpolation of the higher-order derivatives is a natural next step. In this paper, we did not consider an optimality measure for choosing the projection spaces. One might consider combining our framework with the
recent work on optimal parametric model reduction in a joint $\mathcal{H}_2 \otimes \mathcal{L}_2$ measure  \cite{HMMS21}. Even though we have considered here the projection-based approaches, data-driven methods have been also considered for parametric systems \cite{IA14}. Interpreting our reduced model in that framework could provide further hints for data-driven modeling. 

As we stated in \Cref{subsec:numericalConstructionOfVW}, when the power series expansions are truncated, we can no longer guarantee exact interpolation. We will investigate in a future work how the perturbation results from interpolatory model reduction with inexact solves \cite{beattie2012inexact} can be used to quantify the interpolation error due to the truncation.

%\BUcomment{Future work: $\calH_2$-optimal and IRKA, or Wilson conditions), \ldots}

%\BUcomment{Future work II: want to consider different parametrization of $V(p)$ and $W(p)$. Main motivation. If $A(p)\in\R[p]^{\stateDim\times\stateDim}$, and $B(\param)\in\R[p]^{\stateDim\times\inpVarDim}$, then 
%\begin{displaymath}
%	(\lambda I_{\stateDim} - A(p))^{-1}B(p)\in \R(p)^{\stateDim\times 1},
%\end{displaymath}
%i.e., it seams reasonable to search for $V(p),W(p) \in \R(p)^{\stateDim\times\stateDimRed}$.}

%%%%%%%%%%%%%%%%%%%%%%%%%%%%%%%%%%%%%%%%%%%%%%%%%%%%%%%%%%%%%%%%%%%%%%%%%%%%%%%%
\section*{Acknowledgments}
The work of S.~Gugercin
was supported in parts by National Science Foundation under Grant No. 
DMS-1923221 and DMS-1819110.
The work of B.~Unger is funded by the German Research Foundation (DFG) under Germany's Excellence Strategy - EXC 2075 – 390740016. In addition, B.~Unger acknowledges support by the Stuttgart Center for Simulation Science (SimTech). 

\addtolength{\textheight}{-3cm}   % This command serves to balance the column lengths

%%%%%%%%%%%%%%%%%%%%%%%%%%%%%%%%%%%%%%%%%%%%%%%%%%%%%%%%%%%%%%%%%%%%%%%%%%%%%%%%

\bibliographystyle{plain}
\bibliography{literature}    %%%%%%%%%%%%%%%%%%%%%%%%%%%%%%%%%%%%%%%%%%%%%%%%%%%%%%%%%%

\end{document}